\documentclass[12pt]{article}
\usepackage{amssymb,amsmath,amscd,amsthm}

\newtheorem{theorem}{Theorem}[section]

\newtheorem{lemma}[theorem]{Lemma}

\date{}

\begin{document}

\title{Global limit theorems on the convergence of multidimensional random walks to stable processes. }

\author{A. Agbor \footnote{Dept of Mathematics and Statistics, UNC at Charlotte, NC 28223, aagbor@yorktech.edu.
}, S. Molchanov \footnote{Dept of Mathematics and Statistics, UNC at Charlotte, NC 28223, smolchan@uncc.edu.
}, B. Vainberg \footnote{Dept of Mathematics and Statistics, UNC at Charlotte, NC 28223, brvainbe@uncc.edu;
corresponding author.}}

\maketitle

\begin{abstract} Symmetric heavily tailed random walks on $Z^d,~d\geq 1,$ are considered. Under appropriate regularity conditions on the tails of the jump distributions, global (i.e., uniform in $x,t,~|x|+t\to\infty,$) asymptotic behavior of the transition probability $p(t,0,x)$ is obtained. The examples indicate that the regularity conditions are essential.
\end{abstract}
{\bf Key words:} random walk, heavy tail, stable law, large deviations, global asymptotics



 \section {Introduction}
The paper is motivated by two applied problems from the theory of homopolymers and from branching diffusion processes. Both problems can be reduced to asymptotic analysis of solutions to parabolic Schr\"{o}dinger type equation
\begin{equation}\label{Cp}
u_t=-Hu:={\mathcal L}u+V(x)u, \quad u(0,x)=u_0(x).
\end {equation}
Here $t\geq 0,~ x\in Z^d$ (other phase spaces $X$ can be considered, but we will restrict ourselves to the lattice case in this paper), $\mathcal L=\mathcal L^*$ is a self-adjoint  non-positive operator on $L^2(Z^d)$ generating the underlying random walk, and $V$ is a potential.

In the case of homopolymers, problem (\ref{Cp}) with a fast decaying (or even compactly supported) potential $V(x)\geq 0$ describes the shape of a long polymer chain in the field of an attracting potential (see \cite{CrM}, \cite{CrkoMoV}, \cite{CrkoMoV2} and references there).

In the case of branching diffusion processes, $V(x)$ is proportional to the rate of birth of new particles. There are two simplest possibilities: $V(x)=\beta$ is a constant or $V(x)\geq 0$ is compactly supported. In the first case, the rate does not depend on $x$. In particular, this problem appears in the classical model of population introduced (when $X=R^d$) by Kolmogorov, Petrovskii and Piskunov. In the second case, the new particles are produced in a bounded region, see \cite{yar}, \cite{Myar}.

In all the publications mentioned above, as well as in many others, the generator $\mathcal L$ is local. Usually, it is the Laplacian (when $X=R^d$) or lattice Laplacian. However, some more profound models of population dynamics were introduced recently in mathematical biology where the underlying brownian motion in $R^d$ (or simple symmetric   random walk on $Z^d$) is replaced by a heavy- tailed Markov process (containing long jumps, the so-called Levy flights). The main models of this type are contact processes introduced in \cite{kosk}, \cite{kokupi}.
The main goal of the present paper is a study of the global asymptotic behavior of the transition probability for the random walk (time is continuous) with heavy tails. We consider below only lattice models, but it is not difficult to carry the results over to the continuous case.

Let us describe the model under consideration. Let $x(t)$ be the random walk on $Z^d$ with continuous time and the generator
\[
\mathcal Lf(x)=\kappa\sum_{z\in Z^d}(f(x+z)-f(x))a(z), \quad x\in Z^d,
\]
where $a(z)=a(-z)\geq 0$ (we assume that $a$ is symmetric), $ \sum_{z\in Z^d}a(z)=1, ~ \kappa>0$. The process spends an exponentially distributed time $\tau_x$ at each site $x\in Z^d$ ($P\{\tau_x>s\}=\exp(-\kappa s)$) and jumps at time $\tau_x+0$ from $x$ to $x+z$ with the probability $a(z)$. The transition probabilities $p(t,x)=P_x\{x(t)=x;~x(0)=0\}$ satisfy the equation
\[
p_t=\mathcal Lp,\quad p(0,x)=\delta(x).
\]
After the Fourier transform, it takes the form $\widehat{p}_t(k)=\widehat{a}(k)\widehat{p}(k), ~~\widehat{p}(0,k)=1,~$ where $\widehat{a}(k)$ is a periodic function given by
 \begin{equation}\label{ahat}
\widehat{a}(k)=\sum_{z\in Z^d}e^{-i(k,z)}a(z)=\sum_{z\in Z^d}\cos(k,z)a(z).
\end{equation}
Thus
\[
p(t,x)=\frac{1}{(2\pi)^d}\int_{T^d}e^{i(k,x)+\kappa t(\widehat{a}(k)-1)}dk, \quad T^d=[-\pi,\pi]^d.
\]

If $\mathcal L$ is the Laplacian, then $a(z)=1/2d$ when $|z|=1$ and $a(z)=0$ when $|z|\neq1$. In our case, $a(z)$ does not have a compact support, and moreover, we assume that it decays slowly at infinity.

If $E|x(t)|^2<\infty$, i.e., $\sum_{z\in Z^d}|z|^2a(z)<\infty$, then the local ($|x|=O(t^{1/2})$) limit theorem holds \cite{IbL} (under minimal additional assumptions on $x(t)$):
\[
P\{x(t)=x\}\sim \frac{e^{(B^{-1}x,x)/2t\kappa^2}}{(2\pi t\kappa^2)^{d/2}\sqrt{{\rm det}B}}, \quad x\in Z^d, \quad |x|\leq Ct^{1/2},\quad t\to\infty,
\]
where matrix $B$ is the covariance of the jump distribution $a(z)$:
\[
B=[-\frac{\partial^2\widehat{a}(k)}{\partial k_i\partial k_j}]|_{k=0}.
\]
This local limit theorems (in the central zone) were extended \cite{IbL}, \cite{den} for a wide class of stationary and ergodic processes with fast decaying correlations.

We assume that the second moment does not exists, but function $a(z)$ has a regular behavior at infinity:
 \begin{equation}\label{a33}
a(z)=\frac{a_0(\dot{z})}{|z|^{d+\alpha}}(1+o(1)), \quad |z|\to\infty, \quad \dot{z}=\frac{z}{|z|}, ~~a_0(\dot{z})>\delta>0,~~~0<\alpha<2.
\end{equation}
From here it follows that
\begin{equation}\label{c1}
\widehat{a}(k)-1=-b_0(\dot{k})|k|^{\alpha}(1+o(1)), \quad |k|\to 0, \quad \dot{k}=\frac{k}{|k|},
\end{equation}
where $b_0(\dot{k})$ is an appropriate integral transformation of $a_0(\dot{z})$. Relation (\ref{c1}) can be proved under a weaker assumption than (\ref{a33}).

We assumed that $a(-z)=a(z)$ (i.e, the random walk is symmetric). In this case,
\begin{equation}\label{ab1}
b_0(\dot{k})=-\Gamma(-\alpha)\cos\frac{\alpha\pi}{2}\int_{S^{d-1}}a_0(\dot{x})|(\dot{x},\dot{k})|^{\alpha}dS_{\dot{x}}>0,
\end{equation}
where $\Gamma$ is the gamma-function.
Then a different local limit theorem holds:
\begin{equation}\label{6}
P\{x(t)=x\}\sim \frac{1}{t^{d/\alpha}}S(\frac{x}{t^{1/\alpha}}), \quad x\in Z^d, \quad |x|\leq Ct^{1/\alpha},\quad t\to\infty,
\end{equation}
where $S(y)=S_{\alpha,a_0}(y)$ is the density of the stable $d$-dimensional law depending on $\alpha\in (0,2)$ and function $a_0$ in (\ref{a33}). Function $S$ is given by its Fourier transform (characteristic function):
\[
\widehat{S}(k)=e^{-b_0(\dot{k})|k|^\alpha}, \quad k\in R^d.
\]
The local limit theorem (\ref{6}) in the case of $d=1$ can be found in \cite{IbL}, and it can be proved similarly for arbitrary $d$.

Recall (see \cite[Ch XII, \S 11]{Fel}) that the most general symmetric stable law has the characteristic function
\[
\psi_{\alpha,\mu}(k)=e^{-|k|^\alpha\int_{S^{d-1}}|(\dot{x},\dot{k})|^\alpha \mu(d\dot{x})},
\]
where $\mu(d\dot{x})$ is a symmetric (with respect to the reflection over the origin) finite measure on $S^{d-1}$. In our case, the measure has the density $ca_0(\dot{x})\geq \delta>0$.

The goal of this paper is to establish the global asymptotic behavior of $p(t,x)=P\{x(t)=x\}$ when $t+|x|\to\infty$ without restrictions on the relations between $t$ and $|x|$. Thus the zone of large deviations ($|x|=O(t^{1/\alpha+\delta})$ or even $|x|=O(e^{\delta t})$) is also included. One can't expect essential results in  this direction without strong additional assumptions on the tails of $a(z)$ that lead to a more specific behavior of $\widehat{a}(k)$ as $k\to 0$.

 Namely we assume that
\begin{equation}\label{case1}
a(z)= \sum_{j=0}^{d+\epsilon}\frac{a_j(\dot{z})}{|z|^{d+\alpha+j}}+O(\frac{1}{|z|^{2d+\alpha+1+\epsilon}}), \quad |z|\to\infty, ~~~ \alpha\in (0,2), \quad  a_j\in C^{d+1-j+\epsilon}(S^{d-1}),
\end{equation}
where $\epsilon=1$ if $\alpha=1$ and $\epsilon=0$ otherwise, and $a_0(\dot{z})>\delta>0$.

The following lemma (which will be proved in Appendix) provides the asymptotic behavior of $\widehat{a}(k)$ at zero.
\begin{lemma}\label{frt}
If (\ref{case1}) holds,
then
\begin{equation}\label{l01}
\widehat{a}(k)= 1-\sum_{j=0}^{d} b_j(\dot{k})|k|^{\alpha+j}+f(k), \quad k\in T^d=[-\pi,\pi]^d , \quad f(0)=0,
\end{equation}
where $b_j\in C^{d+[\alpha]+1}(S^{d-1})$ and function $f$, being extended periodically on $R^d$, belongs to $ C^{d+[\alpha]+1}(R^d)$.

Moreover, the homogeneous function   $-b_0(\dot{k})|k|^{\alpha}$ in $R^d$ is the Fourier transform of the homogeneous (of order $-d-\alpha$) distribution that is equal to $a_0(\dot{x})|x|^{-d-\alpha}$ when $0\neq x\in R^d$, and
\begin{equation}\label{ab1}
b_0(\dot{k})=-\Gamma(-\alpha)\cos\frac{\alpha\pi}{2}\int_{S^{d-1}}a_0(\dot{x})|(\dot{x},\dot{k})|^{\alpha}dS_{\dot{x}}>0,
\end{equation}
where $\Gamma$ is the gamma-function.
\end{lemma}
{\bf Remarks.} 1) Note that $f$ can not be omitted, since a change in the values of $a(z)$ at several points does not perturb its asymptotic behavior at infinity, but changes $\widehat{a}(k)$ by an analytic function.


2) The next two properties of $\widehat{a}(k)$ follow immediately from properties of $a_0(\dot{x})$:
\begin{equation}\label{cA}
 \widehat{a}(-k)=\widehat{a}(k);\quad  ~\widehat{a}(k)<1, ~~ 0\neq k\in T^d.
\end{equation}
The latter inequality follows from (\ref{ahat}) and the assumption  $\sum_{x\in Z^d}a(x)=1$ if, for each $k\in T^d, ~k\neq 0,$ there is a point $z\in Z^d$, where $e^{-i(x,k)}\neq 1$ and $a(x)\neq 0$. Such points $z$ exist due to (\ref{case1}).

The following uniform asymptotics of the function $p(t,x)=p(t,x,0)$ is one of the main results of the present paper. We put $\kappa=1$ everywhere below, since this condition can be satisfied after simple substitution $\kappa t\to t$.

\begin{theorem}\label{mt}
(1) Let (\ref{l01})-(\ref{cA}) hold. Then
\[
p(t,x)=\frac{1}{t^{d/\alpha}}S(\frac{x}{t^{1/\alpha}})(1+o(1)), \quad {\text when} \quad x\in Z^d, \quad |x|+t\to\infty,
\]
where
\begin{equation}\label{pF}
S(y)=\frac{1}{(2\pi)^{d}}\int_{R^d} e^{i(k,y)-b_0(\dot{k})|k|^{\alpha}}dk>0,~~y\in Z^d,
\end{equation}
is the stable density $S=S_{\alpha,a_0}(y)$.

(2) If $\frac{|x|}{t^{1/\alpha}}\to\infty, ~|x|\geq 1,$ then the previous statement can be specified as follows:
\[
p(t,x)=\frac{a_0(\dot{x})}{t^{d/\alpha}}(\frac{t^{1/\alpha}}{|x|})^{d+\alpha}(1+o(1))=\frac{a_0(\dot{x})t}{|x|^{d+\alpha}}(1+o(1)).
\]
\end{theorem}
{\bf Remark 1.} The study of global limit theorems for the sums $s_n=x_1+...+x_n$ of i.i.d.r.v.  was initiated by Yu. Linnik (see \cite{IbL}, Ch XIV). In 1-D case, he proved a uniform in $x\in R^1$ local limit theorem for the density $p_n(\cdot)$ of $s_n$ under regularity conditions on the tails similar to (\ref{case1}). Some extensions of these results (based on the technique of quasi-cumulants) can be found in \cite{MPS}, \cite{S}.

{\bf Remark 2.} Our proof allows one to write the next terms of the global asymptotics of $p(t,x)$ under additional assumptions on $a(x)$ or $\widehat{a}(k)$ (more terms in the asymptotics (\ref{case1}), (\ref{l01}) and additional smoothness).

The regularity condition (\ref{case1}) on the tails of $a(x)$ usually holds in applications. Therefore, (\ref{l01})-(\ref{cA}) also hold. One has to be careful and check the validity of the latter relations if  (\ref{case1}) is violated.
The following examples show that the smoothness assumptions in (\ref{l01}) are essential.

{\bf Example 1.} The first example concerns a random variable in $R^3$, but it also will be used to construct a specific random walk on $Z^3$ in Example 2.

Consider a random variable with an isotropic distribution density $\alpha (r)=\frac{\sin^4r}{\pi^2r^4}$ in $R^3,~r=|x|.$ Then
\begin{equation}\label{spF}
\widehat{\alpha}(k)=\int_0^\infty\int_{S^2} r^2e^{-ikx}\alpha(r)dSdr=\frac{4\pi}{|k|}\int_0^\infty r\sin( |k|r )\alpha(r)dr,
\end{equation}
where $S^2$ is the unit sphere in $R^3$. We put $\sin^4r=\frac{3}{8}-\frac{1}{2}\cos(2r)+\frac{1}{8}\cos(4r)$ in the expression for $\alpha(r)$ and evaluate the integral using the following classical identity:
\[
\int_0^\infty \frac{\sin (ar)}{r}dr=\frac{\pi}{2}{\rm sgn} a.
\]
This leads to
\[
\widehat{\alpha}(k)=
\begin{cases}
                        1-\frac{3|k|}{8},~~\quad  |k|\leq 2 \\
                        \frac{2}{|k|}-1+\frac{|k|}{8},~~|k|\in[2,4] \\
                        0,~~\quad \quad \quad  |k|>4 .
                      \end{cases}
\]
Function $\widehat{\alpha}(k)$ has the form (\ref{l01}) with $\alpha=1, b_0(\dot{k})=3/8$, and $b_j(\dot{k})=0, j>0,$ but function $f(k)$ for $\widehat{\alpha}(k)$ in our example is not as smooth as Theorem \ref{mt} requires: it belongs to $C^1(R^3)$, but already does not belong to $C^2(R^3)$.

One can evaluate the transition density $p(t,x)=\frac{1}{(2\pi)^3}\int_{R^3}e^{t(\widehat{\alpha}(k)-1)+ixk}dk$ by using spherical coordinates, evaluating the integral over the sphere $S^2$ as in (\ref{spF}), followed by integration by parts three times in $|k|$. This leads to
\[
p(t,x)=\frac{t}{\pi^2r^4}[\frac{3}{8}-\frac{1}{2}\cos(2r)e^{-\frac{3}{4}t}+\frac{1}{8}\cos(4r)e^{-t}]+O(\frac{t}{r^5}),   \quad r\gg t.
\]
This oscillatory behaviour of $p$ is different from the one stated in Theorem \ref{mt}.

{\bf Example 2.} Using the function $\widehat{\alpha}(k)$ defined in the previous example and the Poisson summation formula, one can construct the characteristic function $\widehat{a}(k)$ of the lattice distribution:
\[
\widehat{a}(k)=\sum_{n\in Z^d}\widehat{\alpha}(k+2\pi n).
\]
Since $\widehat{\alpha}(k)=0$ for $|k|>4$, we have $\widehat{a}(k)=\widehat{\alpha}(k)$ on $T^d=[-\pi,\pi]^d$. Thus the jump distribution $a(x)$ and the transition probabilities $p(t,x)$ coincide with the restrictions on the lattice $Z^d$ of the corresponding expressions in Example 1.

{\bf Example 3.} Consider the following well known example: the stable distribution given by the characteristic function
\[
\widehat{S}(k)=e^{-\sum_{j=1}^d |k_j|^\alpha}.
\]
Then
\[
\widehat{a}(k)=1-|k|^\alpha\sum_{j=1}^d |\dot{k}_j|^\alpha+O(|k|^{1+\alpha}),  \quad k\to 0,
\]
i.e., function $b_0(\dot{k})$ is not smooth. Since $S(x)=\prod_{j=1}^d S_\alpha(x_j)$, the asymptotics of $S(x)$ at infinity depends significantly on $\dot{x}$. Say,
\[
S(r,0,...,0)\sim \frac{c}{r^{1+\alpha}}, \quad r\to\infty,
\]
but
\[
S(r,r,...,r)\sim \frac{c^d}{r^{d(1+\alpha)}}, \quad r\to\infty.
\]

Note that the problem of the asymptotic behaviour for general stable distributions is not simple. There are many relatively recent publications on that topic, see \cite{ark}, \cite{hi}, \cite{wat} and references there.
\section{Proof of the main result}
The proof will use three lemmas proved below and the following simple facts. Since function $\widehat{a}(k)$ is even, we have
\begin{equation}\label{grf}
\nabla f(k)=0 \quad {\rm at} \quad k=0.
\end{equation}
From
 (\ref{c1}), (\ref{ab1}) and (\ref{cA}) it follows that
\begin{equation}\label{c0}
\widehat{a}(k)-1<-\gamma|k|^\alpha,~~~k\in[-\pi,\pi]^d,~~~\gamma>0.
\end{equation}

Let $\psi=\psi(\tau)\in C^\infty(R^1)$ be a cut-off function such that $\psi(\tau)=1,~|\tau|<1,~~\psi(\tau)=0,~|\tau|>2$. Let
\[
I=\int_{[-\pi,\pi]^d}e^{i(k,x)+t(\widehat{a}(k)-1)}[1-\psi(|k|t^{1/\alpha})]dk.
\]
\begin{lemma}\label{bes}
Let (\ref{case1}) hold. Then the following estimate is valid for all $t, |x|>0$:
\[
|I|\leq \frac{C}{t^{d/\alpha}}(\frac{t^{1/\alpha}}{|x|})^{d+[\alpha]+1}.
\]
\end{lemma}
{\bf Proof.} From (\ref{c0}), (\ref{c1}) it follows that
\begin{equation}\label{dife}
|\partial^j_ke^{t(\widehat{a}(k)-1)}|\leq C\sum_{l=1}^{|j|}(|k|^{l\alpha-|j|}t^l)e^{-\gamma|k|^\alpha t}, \quad k\in[-\pi,\pi]^d , ~~0<|j|\leq d+[\alpha]+1.
\end{equation}
Indeed, $\frac{\partial}{\partial k_s}e^{t(\widehat{a}(k)-1)}=\frac{\partial \widehat{a}(k)}{\partial k_s}te^{t(\widehat{a}(k)-1)}$, and the pre-exponential factor can be estimated from above by $C|k|^{\alpha-1}t$, due to (\ref{c1}), (\ref{grf}). Each next differentiation leads to the appearance of an additional similar  pre-exponential factor (when the derivative is applied to the exponent) or to decreasing of the power of $k$ in the estimate of the pre-exponent by one. These arguments justify (\ref{dife}).

Since $|k|^{l\alpha} t^le^{-\gamma|k|^\alpha t/2}\leq C$, estimate (\ref{dife}) implies that
\[
|\partial^j_ke^{t(\widehat{a}(k)-1)}|\leq Ct^{|j|/\alpha}e^{-\gamma |k|^\alpha t/2}, \quad |k|t^{1/\alpha}>1 , ~~|j|\leq d+[\alpha]+1.
\]
Hence, for each $s=1,2,..,d$ and $m=d+[\alpha]+1$, we have
\[
|\frac{\partial^m}{(\partial k_s)^m}[e^{t(\widehat{a}(k)-1)}(1-\psi(|k|t^{1/\alpha}))]|\leq C t^{m/\alpha} e^{-\gamma|k|^\alpha t/2},
\]
and therefore (since $\widehat{a}(k)$ is periodic),
\[
|I|=|(-ix_s)^{-m}\int_{[-\pi,\pi]^d}e^{i(k,x)}\frac{\partial^m}{(\partial k_s)^m}[e^{t(\widehat{a}(k)-1)}(1-\psi(|k|t^{1/\alpha}))]dk|
\]
\[
\leq \frac{Ct^{m/\alpha}}{|x_s|^m}\int_{R^d}e^{-\gamma |k|^\alpha t/2}dk=\frac{C_1t^{(m-d)/\alpha}}{|x_s|^m}.
\]
This completes the proof of the lemma, since $s$ is arbitrary.

\qed

Let
\[
I_1=\frac{1}{(2\pi)^{d}}\int_{|k|t^{1/\alpha}<2}e^{i(k,x)+t(\widehat{a}(k)-1)}\psi(|k|t^{1/\alpha})dk.
\]
\begin{lemma}\label{bes1}
Let (\ref{case1}) hold. Then for each $r>0, ~|x|>rt^{1/\alpha}$ and $|x|+t\to \infty$, the following asymptotics holds when $\frac{|x|}{t^{1/\alpha}}\to\infty$
\[
I_1= \frac{a_0(\dot{x})}{t^{d/\alpha}}(\frac{t^{1/\alpha}}{|x|})^{d+\alpha}(1+o(1))=\frac{a_0(\dot{x})t}{|x|^{d+\alpha}}(1+o(1)),
\]
where there exist $\varepsilon>0$ and for each $r>0$ there exists $C=C(r)$ such that
\[
|o(1)|<C(|x|^{-\varepsilon}+(\frac{t^{1/\alpha}}{|x|})^{\varepsilon}),~~ ~~|x|>rt^{1/\alpha},~~|x|>1.
\]
\end{lemma}
{\bf Proof.} The substitution $kt^{1/\alpha}\to k$ followed by integration by parts implies that, for each $N>0$,
\[
|\int_{|k|t^{1/\alpha}<2}e^{i(k,x)}\psi(|k|t^{1/\alpha})dk|<C_N\frac{1}{t^{d/\alpha}}(\frac{t^{1/\alpha}}{|x|})^N=C_N\frac{t}{|x|^{d+\alpha}}(\frac{t^{1/\alpha}}{|x|})^{N-d-\alpha}
=\frac{t}{|x|^{d+\alpha}}o(1).
\]
Thus, it is enough to proof Lemma \ref{bes1} for
\[
I_2=\frac{1}{(2\pi)^{d}}\int_{|k|t^{1/\alpha}<2}e^{i(k,x)}(e^{t(\widehat{a}(k)-1)}-1)\psi(|k|t^{1/\alpha})dk.
\]

From  (\ref{case1}) it follows that
\begin{equation}\label{del}
t(\widehat{a}(k)-1)=-tb_0(\dot{k})|k|^{\alpha}+O(t|k|^{\alpha+\delta}),
\end{equation}
where $\delta=1 $ if $\alpha\leq 1$ and $\delta=2-\alpha $ if $\alpha> 1$. Each term on the right is bounded when $|k|t^{1/\alpha}<2$. Thus
\begin{equation}\label{expp}
e^{t(\widehat{a}(k)-1)}-1=-tb_0(\dot{k})|k|^{\alpha}+O(t|k|^{\alpha+\delta})+O(t^2|k|^{2\alpha}), \quad |k|t^{1/\alpha}<2.
\end{equation}
We replace the middle factor in the integrand of $I_2$ by the right-hand side in (\ref{expp}) and split the integral into three terms. Then we estimate the last two integrals using the substitution $k\to |x|^{-1}k$. The domain of integration remains bounded after the substitution (since $|x|>rt^{1/\alpha}$), and therefore the last two integrals can be estimated by
\[
C(t|x|^{-d-\alpha-\delta}+t^2|x|^{-d-2\alpha})=\frac{t}{|x|^{d+\alpha}}o(1).
\]
Now it remains to prove the statement of the lemma for the integral
\[
I_3=\frac{-t}{(2\pi)^{d}}\int_{|k|t^{1/\alpha}<2}e^{i(k,x)}b_0(\dot{k})|k|^{\alpha}\psi(|k|t^{1/\alpha})dk=t\mathcal F^{-1} (-b_0(\dot{k})|k|^{\alpha}\psi(|k|t^{1/\alpha})).
\]
where $\mathcal F=\mathcal F_{x\to k}$ is the Fourier transform operator.

Due to Lemma \ref{frt}, $t\mathcal F^{-1} (-b_0(\dot{k})|k| ^{\alpha})=a_0(\dot{x})t|x|^{-d-\alpha}$. Thus it is enough to show that
\[
\mathcal F^{-1} [-b_0(\dot{k})|k|^{\alpha}(1-\psi(|k|t^{1/\alpha}))]=|x|^{-d-\alpha}o(1).
\]
The latter estimate can be obtain by the substitution $k\to kt^{-1/\alpha}$ in the left-hand side followed by integration by parts $d+[\alpha]+1$ times with respect to each of the variables $k_l,~1\leq l\leq d$. This leads to the estimate of the left-hand side by
\[
Ct^{\frac{-d-\alpha}{\alpha}}(\frac{t^{1/\alpha}}{|x|})^{d+[\alpha]+1}=\frac{C}{|x|^{d+\alpha}}(\frac{t^{1/\alpha}}{|x|})^{[\alpha]+1-\alpha}=|x|^{-d-\alpha}o(1).
\]
\qed

\begin{lemma}\label{l3}
Function $S$ defined in (\ref{pF}) is positive.
\end{lemma}
{\bf Proof.} Consider first the case of $\alpha\in(1,2)$. It will be shown below (see (\ref{pp0})) that
\[
S(y)=\lim_{t\to\infty} t^{d/\alpha}p(t,t^{1/\alpha}y).
\]
Thus $S(y)\geq 0$. From (\ref{pF}) it follows that $S(0)>0. $ Let $|y|<R$ be the largest ball where $S(y)>0$ for all $y$ inside the ball. We need to show that $R=\infty$.

Consider the convolution $S_2=S\ast S.$ Its Fourier transform is equal to $e^{-2b_0(\dot{k})|k|^{\alpha}}$, and therefore
\begin{equation}\label{pos}
S_2(y)=\frac{1}{(2\pi)^{d}}\int_{R^d}e^{i(k,x)-2b_0(\dot{k})|k|^{\alpha}}dk=2^{-\frac{d}{\alpha}}S(2^{\frac{1}{\alpha}}y).
\end{equation}
The last equality (self-similarity) can be obtained by substitution $k\to 2^{1/\alpha}k$. Since  $S(y)\geq 0$, the convolution is positive in the ball $|y|<2R$. Hence (\ref{pos}) implies that $S(y)> 0$ in the ball $|y|<2^{1-\frac{1}{\alpha}}R$. The radius of the latter ball is larger than $R$ if $\alpha\in (1,2),~R<\infty$. Thus $R=\infty$, i.e., $S(y)>0$ everywhere if $\alpha\in (1,2)$.

Let now $0<\alpha\leq 1$. By passing in (\ref{pF}) to spherical coordinates and integrating by parts, we obtain that
\[
S(y)=\frac{1}{(2\pi)^{d}}\int_{S^{d-1}}\int_0^\infty \cos[|k||y||(\dot{k},\dot{y})|]e^{-b_0(\dot{k})|k|^{\alpha}}d|k|dS_{\dot{k}}
\]
\[
=\frac{1}{(2\pi)^{d}}\int_{S^{d-1}}\int_0^\infty \frac{\sin[|k||y||(\dot{k},\dot{y})|]}{|y|(\dot{k},\dot{y})}b_0(\dot{k})\alpha|k|^{\alpha-1}e^{-b_0(\dot{k})|k|^{\alpha}}d|k|dS_{\dot{k}}.
\]
The last integral is positive since the function $b_0|k|^{\alpha-1}e^{-b_0(\dot{k})|k|^{\alpha}}$ is monotone in $|k|$, and the integral
\[
\int_0^\infty\sin(\gamma\tau)f(\tau)d\tau
\]
is positive for every $\gamma>0$ and every monotonically decreasing function $f$.

\qed

{\bf Proof of Theorem \ref{mt}.} The second statement of the theorem follows immediately from Lemmas \ref{bes}, \ref{bes1}. It remains to prove the first statement.

First, let us prove that the second statement of the theorem is valid for the function
\[
p_0(t,x)=\frac{1}{(2\pi)^{d}}\int_{R^d} e^{i(k,x)-tb_0(\dot{k})|k|^{\alpha}}dk=\frac{1}{t^{d/\alpha}}S(\frac{x}{t^{1/\alpha}}).
\]
We split $p_0$ in two parts $p_0(t,x)=p_{0,1}+p_{0,2}$, where
\[
p_{0,1}=\frac{1}{(2\pi)^{d}}\int_{R^d} e^{i(k,x)-b_0(\dot{k})|k|^{\alpha}}\psi(kt^{1/\alpha})dk,~~~p_{0,2}=\frac{1}{(2\pi)^{d}}\int_{R^d} e^{i(k,x)-b_0(\dot{k})|k|^{\alpha}}(1-\psi(kt^{1/\alpha}))dk,
\]
The estimate for $I$, established in Lemma \ref{bes}, is valid for $p_{0,2}$. We can not formally refer to Lemma \ref{bes} (since the integrals are different), but the same arguments can be applied. In fact, they are simpler in the case of  $p_{0,2}$. We change the variable $k\to kt^{-1/\alpha}$ and then integrate by parts $d+[\alpha]+1 $ times with respect to each of the variables $k_l,~1\leq l\leq d$. This leads to
\[
|p_{0,2}|\leq \frac{C}{t^{d/\alpha}}(\frac{t^{1/\alpha}}{|x|})^{d+[\alpha]+1}=C\frac{t}{|x|^{d+\alpha}}(\frac{t^{1/\alpha}}{|x|})^{[\alpha]+1-\alpha}=\frac{t}{|x|^{d+\alpha}}o(1),
\]
where $|o(1)|\to 0$ as $\frac{|x|}{t^{1/\alpha}}\to\infty$. Lemma \ref{bes1} is valid when $\widehat{a}(k)-1=-b_0(\dot{k})|k|^{\alpha}$ and can be applied to $p_{0,1}$. Together with the estimate for  $p_{0,2}$ this justifies that the second statement of the theorem is valid for $p_0$. Thus $p_0$ and $p$ are close to the same non-zero function $\frac{a(\dot{x})t}{|x|^{d+\alpha}}$, and therefore are close to each other. Hence, for each $\varepsilon>0$, there exists an $A=A(\varepsilon)$ such that
\[
|p(t,x)-p_0(t,x)|\leq \varepsilon p_0(t,x) \quad{\rm when} \quad \frac{|x|}{t^{1/\alpha}}\geq A, ~|x|\geq A.
\]
Now the theorem will be proved if we show that for each fixed $\varepsilon, A>0$ there exists  $T=T(\varepsilon, A)$ such that
\begin{equation}\label{lart}
|p(t,x)-p_0(t,x)|\leq \varepsilon p_0(t,x) \quad{\rm when} \quad \frac{|x|}{t^{1/\alpha}}\leq A, ~t>T.
\end{equation}

We will show that
\begin{equation}\label{pp0}
|p(t,x)-p_0(t,x)|\leq Ct^{-d/\alpha-\varepsilon}, ~~\varepsilon>0,     \quad{\rm when} \quad \frac{|x|}{t^{1/\alpha}}\leq A,~ ~t\to\infty.
\end{equation}
Then we note that $p_0$ has the form $p_0(t,x)=t^{-d/\alpha}S(\frac{|x|}{t^{1/\alpha}})$, where $S$ is given by (\ref{pF}). Since the function $b_0(\dot{k})$ is even, $S$ is real-valued. Since $S=S(y)$ is continuous in $y$ and does not vanish (see Lemma \ref{l3}
), $p_0>c_0t^{-d/\alpha}>0$ when $|y|\leq A$. Thus (\ref{pp0}) implies (\ref{lart}), so it is enough to justify (\ref{pp0}).

Let
 \[
p^{(1)}(t,x)=\frac{1}{(2\pi)^{d}}\int_{|k|<t^{-1/(\alpha+\delta)}} e^{i(k,x)+t(\widehat{a}(k)-1)}dk,
 \]
 \[
 p^{(1)}_0(t,x)=\frac{1}{(2\pi)^{d}}\int_{|k|<t^{-1/(\alpha+\delta)}} e^{i(k,x)-tb_0(\dot{k})|k|^\alpha}dk,
\]
where $\delta$ was defined in (\ref{del}). Let $p^{(1)}(t,x), ~ p^{(1)}_0(t,x)$ be the same integrals over the complementary sets in $T^d$ and $R^d$, respectively. From (\ref{c0}) it follows that $t(\widehat{a}(k)-1)<-\gamma t^{\delta/(\alpha+\delta)}$, and therefore $|p^{(2)}|<Ce^{-\gamma t^{\delta/(\alpha+\delta)}}$. Obviously, similar estimate is valid for $p^{(2)}_0$. It remains to prove (\ref{pp0}) for $p^{(1)}(t,x)-p^{(1)}_0(t,x)$.

Due to (\ref{del}), we have
\[
e^{t(\widehat{a}(k)-1)}-e^{-tb_0(\dot{k})|k|^\alpha}=e^{-tb_0(\dot{k})|k|^\alpha}O(t|k|^{\alpha+\delta}), \quad  {\rm when} ~~|k|<t^{-1/(\alpha+\delta)}.
\]
Hence
\[
|p^{(1)}-p^{(1)}_0|<C\int_{|k|<t^{-1/(\alpha+\delta)}}e^{-tb_0(\dot{k})|k|^\alpha}t|k|^{\alpha+\delta}dk
<C\int_{R^d}e^{-tb_0(\dot{k})|k|^\alpha}t|k|^{\alpha+\delta}dk=C_1t^{-\frac{d+\delta}{\alpha}},
\]
and this justifies the validity of (\ref{pp0}) for $p^{(1)}-p^{(1)}_0$.

\qed

\begin{center}
{\bf Appendix}
\end{center}
\begin{lemma}\label{atob}
Let $f\in {\mathcal S}'$ be a homogeneous distribution equal to $u(\dot{x})|x|^{-d-\gamma}$ when $x\neq 0$, where $\gamma>0,~u\in C^{n_1}(S^{d-1})$. Then its Fourier transform $\widetilde{f}={\mathcal F}_{x\to k}f$ is equal to $c(\dot{k})|k|^\gamma$, where $c(\dot{k})\in C^{n_2}(S^{d-1})$ with $n_2=n_1+[\gamma]$ if $\gamma$ is non-integer, and $n_2=n_1+\gamma-1$ if $\gamma$ is an integer.

If $u$ is even (i.e., $u(-\dot{x})=u(\dot{x})$), then
\begin{equation}\label{23a}
c(\dot{k})=\Gamma(-\gamma)\cos\frac{\gamma\pi}{2}\int_{S^{d-1}}u(\dot{x})|(\dot{x},\dot{k})|^{\gamma}dS_{\dot{x}},
\end{equation}
where $\Gamma$ is the gamma-function.
\end{lemma}
 {\bf Proof.} Obviously, $\widetilde{f}$ is homogeneous, i.e., $\widetilde{f}=b(\dot{k})|k|^\gamma$. Let $\varsigma=\varsigma (x)\in C^\infty(R^d),~\varsigma=1 $ for $|x|>2,~\varsigma=0 $ for $|x|<1$, and let $g_l=\frac{\partial^{n_1}}{(\partial x_l)^{n_1}}(\varsigma f)$. Then $|g_l|<C(1+|x|)^{-d-\gamma-n_1}$, and therefore ${\mathcal F}g_l\in C^{n_2}(R^d),~1\leq l\leq d.$ Hence
\[
k_l^{n_1}{\mathcal F}(\varsigma f)\in C^{n_2}(R^d),~~1\leq l\leq d.
\]
This implies that ${\mathcal F}(\varsigma f)|_{|k|=1}\in C^{n_2}(S^{d-1})$, and therefore $b\in C^{n_2}(S^{d-1})$, since ${\mathcal F}((1-\varsigma) f)$ is analytic in $k$.

The proof of (\ref{23a}) can be considered as an exercise in the theory of distributions. In order to prove  (\ref{23a}), we note that $f$ and its Fourier transform are defined for negative $\gamma$. Consider $\gamma\in (-1,0)$. We write the integral ${\mathcal F}f$ in spherical coordinates and evaluate the integral in $|x|$. This leads to (\ref{23a}) when  $\gamma\in (-1,0)$. For other values of $\gamma$, this formula is obtained by analytic continuation in $\gamma$.

\qed

{\bf Proof of Lemma \ref{frt}.} Let $x\in R^d,~z\in Z^d$. We split $R^d$ into unit cubes
$$
q_z:=\{x=z+\tau: 0\leq \tau_s<1, ~1\leq s\leq d\}.
$$
Let $u=u(\dot{x})\in C^n(S^{d-1})$. The following asymptotic expansion holds when $\gamma>0,~x\in q_z,~|z|\to \infty$:
\[
u(\dot{x})|x|^{-d-\gamma}=\sum_{s=0}^{n-1} c_s(\dot{z},\tau)|z|^{-d-\gamma-s}+O(|z|^{-d-\gamma-n}), \quad c_j\in C^{n-j}, \quad c_0=u(\dot{z}).
\]
 Thus
\[
u(\dot{x})|x|^{-d-\gamma}e^{-i(x,k)}=e^{-i(z,k)}[\sum_{s=0}^{n-1} c_s(\dot{z},\tau)|z|^{-d-\gamma-s}e^{-i(\tau,k)}+O(|z|^{-d-\gamma-n})].
\]
We integrate both sides of this relation over the unit cube in $\tau$. The coefficients on the right will be $n-j$ times differentiable functions of $\dot{z}$ and $k$ and they will be  analytic in $k$. We replace variables $k$ by $\sigma=(\sigma_1,...,\sigma_d),~\sigma_l=e^{ik_l}-1,~1\leq l\leq d$. This leads to the following asymptotic expansion
\begin{equation}\label{23}
\int_{q_z}u(\dot{x})|x|^{-d-\gamma}e^{-i(x,k)}dx= e^{-i(z,k)}[\sum_{s=0}^{n-1} d_s(\dot{z},\sigma)|z|^{-d-\gamma-s}+O(|z|^{-d-\gamma-n})], \quad |z|\to \infty,
\end{equation}
where $d_s$ is $n-s$ times differentiable in variables $\dot{z}, \sigma$ and analytic in $\sigma$, and $d_0(\dot{z},0)=u(\dot{z})$.

Next, we write each function $d_s$ in the form $d_s(\dot{z},\sigma)=d_s(\dot{z},0)+\sum_{l=1}^d\sigma_ld_{s,l}(\dot{z},\sigma),$
where $d_{s,l}$ are analytic in $\sigma$. Then we apply similar expansion to functions $d_{s,l}$ and repeat this procedure $n-s$ times in total. This leads to the following representation for $d_s$:
\begin{equation}\label{byp}
d_s(\dot{z},\sigma)=d_s(\dot{z},0)+p_{s,1}(\dot{z},\sigma)+p_{s,2}(\dot{z},\sigma),
\end{equation}
where $p_{s,1}$ are polynomials in $\sigma$ of order $n-s-1$ without the free term and $p_{s,2}$ are homogeneous polynomials in $\sigma$ of order $n-s$ whose coefficients depend on both $\dot{z}$ and $\sigma$. Then we take summation in both sides of (\ref{23}) over all $z\in Z^d$ except some neighborhood of the origin (to avoid singularities at $z=0$) and apply the following relation to each factor $\sigma_l$ in polynomials $p_{s,1},~p_{s,2}$:
\begin{equation}\label{234}
\sum_{Z^d}r(z)e^{-i(z,k)}\sigma_l=\sum_{Z^d}r(z)e^{-i(z,k)}(e^{ik_l}-1)=\sum_{Z^d}(r(z-z^{(l)})-r(z))e^{-i(z,k)},
\end{equation}
where all the components of $z^{(l)}$ are zeroes except $z_l$, which is equal to one. Relation (\ref{234}) shows that the multiplication by $\sigma_l$ reduces every series with the terms of the form $v(\dot{z})|z|^{-d-\gamma'}e^{-i(z,k)}$ (with an arbitrary $\gamma'>0$ and a smooth $v=v(\dot{z})$) to a series whose terms decay faster at infinity. Then from (\ref{23}) it follows that
\begin{equation}\label{rec0}
\int_{|x|>1}u(\dot{x})|x|^{-d-\gamma}e^{-i(x,k)}dx=\sum_{0\neq z\in Z^d} e^{-i(z,k)}[\sum_{s=0}^{n-1}v_s(\dot{z})|z|^{-d-\gamma-s}+O(|z|^{-d-\gamma-n})]+g(k),
\end{equation}
where $v_s\in C^{n-s},~v_0(\dot{z})=u(\dot{z})$, and $g$ is an analytic in $k$ function that appeared because we did not care to evaluate the contribution of the terms with small values of $z$ and the contribution of the integral over the area in a neighborhood of the origin explicitly.

We will use (\ref{rec0}) only for $\gamma$ and $n$ such that
\begin{equation}\label{gn}
\gamma+n=d+\alpha+1 \quad {\rm when }~~\alpha\in (0,2),~~\alpha\neq 1; \quad  \gamma+n=d+3 \quad {\rm when }~~\alpha=1.
\end{equation}
Then the sum in $z$ on the right that contains the remainder term is $m:=d+[\alpha]+1$ times differentiable in $k$, and therefore (\ref{rec0}) implies that
\begin{equation}\label{rec00}
\int_{|x|>1}u(\dot{x})|x|^{-d-\gamma}e^{-i(x,k)}dx=\sum_{0\neq z\in Z^d} e^{-i(z,k)}\sum_{s=0}^{n-1}v_s(\dot{z})|z|^{-d-\gamma-s}+q(k), \quad q\in C^m.
\end{equation}
The left-hand side in (\ref{rec00}) differs from the Fourier transform $c(\dot{k})|k|^{\gamma}$ of $u(\dot{x})|x|^{-d-\gamma}$ (see Lemma \ref{atob}) by an analytic function. Thus (\ref{rec00}) leads to the following expression for the term on the right with $s=0$:
\begin{equation}\label{ind}
\sum_{0\neq z\in Z^d} e^{-i(z,k)}u(\dot{z})|z|^{-d-\gamma}=c(\dot{k})|k|^{\gamma}-\sum_{0\neq z\in Z^d} e^{-i(z,k)}\sum_{s=1}^{n-1}v_s(\dot{z})|z|^{-d-\gamma-s}+h(k), ~~h\in C^m,
\end{equation}
where $c(\dot{k})$ and $u(\dot{z})$ are related by (\ref{23a}). From (\ref{gn}) and Lemma \ref{atob}  (with $n_1=n$) it follows that $c\in C^m$, i.e., the smoothness of functions $c$ and $h$ does not depend on $n$ and $\gamma$ (under condition (\ref{gn})). This allows us to apply the induction on $n=1,2,...$ in formula (\ref{ind}), which leads to the following result:
\begin{equation}\label{ind1}
\sum_{0\neq z\in Z^d} e^{-i(z,k)}u(\dot{z})|z|^{-d-\gamma}=\sum_{s=0}^{n-1}c_s(\dot{k})|k|^{\gamma-s}+f(k), ~~f\in C^m,
\end{equation}
where $c_s\in C^m$, $c_0=c$ is given by (\ref{23a}), and $\gamma$ and $n$ are related by (\ref{gn}). In particular, (\ref{ind1}) holds for $\gamma=\alpha+j, ~u=a_j$ for all the values of $j$ from (\ref{case1}), and the summation of  (\ref{ind1}) in $j$ implies the statement of the Lemma.

\qed


\begin{thebibliography}{102}

\bibitem{ark} S.V. Arhipov, Density function's asymptotic representationin the case of multidimentional strictly stable distribution, In V. Kalashnikov and
V. Zolotarev (Eds.), Lecture Notes in Math. (1987), Volume 1412, pp. 1-21. Springer.

 \bibitem{CrM} M. Cranston, S. Molchanov, On phase transitions and limit theorems for homopolymers. Probability and mathematical physics, CRM Proc. Lecture Notes, 42 (2007), 97–112

\bibitem{CrkoMoV} M. Cranston, L. Koralov, S. Molchanov, B. Vainberg, Continuous model for homopolymers. J. Funct. Anal.  256  (2009),  no. 8, 2656-2696.


\bibitem{CrkoMoV2} M. Cranston, L. Koralov, S. Molchanov, B. Vainberg, A solvable model for homopolymers and self-similarity near the critical point. Random Oper. Stoch. Equ.  18  (2010),  no. 1, 73-95.

 \bibitem{Fel}   W. Feller, Introduction to Probability Theory and its Applications, Vol. 2,  Wiley, (1971).

 \bibitem{hi} S. Hiraba, Asymptotic estimates for densities of muti-dimensional stable distribution, Tsukuba J. Math. 27 (2003), 261-287.

\bibitem{IbL}    I. Ibragimov, Yu. Linnik, Independent and Stationary Sequences of Random Variables, Wolters-Noordhoff publishers, Groningen, Nitherlands, 1971

\bibitem{den}   A. Jakubowski, M. Denker, Stable limit distributions for strongly mixed sequences, Stat. Probab., 8 (1989), 477-483.

\bibitem{kokupi} Yu. Kondratiev, O. Kutoviy, S. Pirogov, Correlation functions and invariant measure
in continuous contact model. Inf. Dimens. Analysis, Quantum probab. And Related
Topics (11)(2008), No 2, 231-258.

\bibitem{kosk} Yu.Kondratiev, A. Skorohod, On contact processes in continium, Inf. Dimensional Analysis, Quantum Probab. and Related Topics, 9, No 2 (2006) , 187-198.

\bibitem{kormol}  L. Koralov, S. Molchanov, Structure of population inside propagating front. Problems in mathematical analysis. No. 69. J. Math. Sci. (N. Y.)  189  (2013),  no. 4, 637-658.

\bibitem{kk}Yu. Kondratiev, O. Kutoviy, S. Molchanov,  On the population dynamics in the stationary random environment, Markov Processes and Related Fields, accepted.

\bibitem{MPS} S. Molchanov, V. Petrov, N. Squartini, Quisiqumulants and limit theorems in the case of Cauchy limiting law, Markov Processes and Related Fields, (2007), 13, 597-624.

\bibitem{yar} S. Molchanov,  E. Yarovaya, The branching processes with lattice space dynimics and finite number of the sets of the generation of the particles, Dokl. Russ. Acad. Sci., Math., 446, No 3 (2012), 259-262.

\bibitem{Myar} S. Molchanov, E. Yarovaya, Structure of the population inside the propagating front of a branching process with a finite set of generating centers, Dokl. Russ. Acad. Sci., Math., 447,  No. 3  (2012), 265--268;  translation in  Dokl. Math.  86  (2012),  no. 3, 787-790.

\bibitem{S} N. Squartini, Global limit theorems for sums of i.i.d.r.v., using quasicummulants, PhD thesis, UNC Charlotte, 2007.

\bibitem{wat} T. Watanaba,  Asymptotic estimates of muti-dimensional stable densities and their applications, Transections of AMS,
359, No 6, (2007), 2851-2879.

\end{thebibliography}
\end{document}